\documentclass[a4paper,12pt]{amsart}
\usepackage{amsmath}
\usepackage{amssymb}
\usepackage{ifthen}
\usepackage{graphicx}
\usepackage{float}
\usepackage{cite}
\usepackage{amsfonts}
\usepackage{amscd}
\usepackage{amsxtra}
\usepackage{color}

\newtheorem{theorem}{Theorem}[section]
\newtheorem{lemma}[theorem]{Lemma}
\newtheorem{corollary}[theorem]{Corollary}

\theoremstyle{definition}

\numberwithin{equation}{section}

\def\be{\begin{equation}}
\def\ee{\end{equation}}

\newcounter{alphabet}


\begin{document}

\title[Initial Successive coefficients of Inverse functions]{Initial Successive coefficients of Inverse functions of certain classes of univalent functions}

\author[Vasudevarao Allu]{Vasudevarao Allu}
\address{School of Basic Science \\ Indian Institute of Technology Bhubaneswar\\
	India}
\email{avrao@iitbbs.ac.in}
\author[Vibhuti Arora]{Vibhuti Arora}
\address{School of Basic Science \\ Indian Institute of Technology Bhubaneswar\\
	India}
\email{vibhutiarora1991@gmail.com}

\subjclass[2010]{30D30, 30C45, 30C50 30C55.}
\keywords{Inverse coefficients, Successive coefficients, Univalent functions}

\begin{abstract}
We consider functions of the type $f(z)=z+a_2z^2+a_3z^3+\cdots$ from a
family of all analytic and univalent functions in the unit disk. Let $F$ be the inverse function of $f$, given by $F(z)=w+\sum_{n=2}^{\infty}A_nw^n$ defined on some $|w|\le r_0(f)$. In this paper, we find the sharp bounds of $\big | |A_{n+1}|-|A_n|\big |$, for $n=1,\,2$, for some subclasses of univalent functions. 
\end{abstract}

\maketitle

\section{Introduction}\label{Introduction}

Let $\mathcal{A}$ denote the class of functions $f$ analytic in the unit disk $\mathbb{D}:=\{z\in \mathbb{C}:|z|<1\}$ with Taylor series
\begin{equation}\label{S}
f(z)=a_1z+a_2z^2+a_3z^3+\cdots,
\end{equation}
with $a_1=1$. Let $\mathcal{S}$ be the set all functions $f\in \mathcal{A}$ that are univalent in $\mathbb{D}$.  Denote by $\mathcal{S}^*$, the family of functions $f$ in $\mathcal{A}$ such that $f(\mathbb{D})$ is a starlike domain with respect to the origin. The family of all functions $f \in \mathcal{A}$ for which $f(\mathbb{D})$ is a convex domain is denoted by $\mathcal{C}$ (see \cite{Dur83,TTV18}).
In 1985, de Branges \cite{DeB1} solved the popular Bieberbach conjecture, which was conjectured in $1916$ by Bieberbach, which states that the Taylor coefficients $a_n$ of functions $f\in \mathcal{S}$ of the form \eqref{S} satisfy the inequality $|a_n| \le n$ and furthermore, equality could
occur if $f$ is some rotation of the Koebe function $k(z):=z/(1-z)^2$. Similarly, the problem
of estimating sharp bound for successive coefficients, namely,
$|a_{n+1}|-|a_n|$, is also an
interesting coefficient problem for a function to be in class $\mathcal{S}$. This problem was first studied
by Goluzin \cite{Gol46} with an idea to solve the Bieberbach conjecture. Hayman \cite{Hay63} proved $\big||a_{n+1}|-|a_n|\big |\le A$ for $f \in \mathcal{S}$, where $A \ge 1$
is an absolute constant and the best known estimate as of now is $3.61$ due to Grinspan \cite{Gri76}. On the other hand, for the class $\mathcal{S}$ sharp bound is known only for $n=2$ (see \cite[Theorem~3.11]{Dur83}), namely
$$
-1\leq|a_3|-|a_2|\leq 1.029\ldots.
$$For convex functions, Li and Sugawa \cite{LS17}  investigated the sharp upper bound of $|a_{n+1}|-|a_n|$ for $n\ge 2$, and sharp lower bounds for $n=2,3$.
Several results are known in this direction \cite{PO19}. These observations are also addressed in the recent papers (see \cite{APS19, Vibhuti}).\\[2mm]

For $f\in \mathcal{S}$ denote by $F$ the inverse of $f$ given by
\begin{equation*}
F(w)=w+\sum_{n=2}^{\infty}A_nw^n,
\end{equation*}
valid on some disk $|w|\le r_0(f)$.
Since $f(f^{-1}(w))=w$, we can easily obtain by equating the coefficients
\begin{equation}\label{A2}
	A_2=-a_2 \text{ and } A_3=2a_2^2-a_3.
\end{equation}
The inverse functions are studied by several authors in different perspective (see, for instance, \cite{TTV18, SimT20} and reference therein).

Although, the sharp bounds are known for $|A_n|$, for $n\ge 2$, when $f\in \mathcal{S}$ (see \cite{Lowner}), but the successive coefficient problem for inverse functions, {\it i.e.,} the bounds $|A_{n+1}|-|A_n|$, is still not known for several important class of functions including the whole class $\mathcal{S}$. So, it seems reasonable and interesting to compute the bounds of $|A_{n+1}|-|A_n|$ for the class of univalent functions and its subclasses even for some particular values of $n$. In \cite{SimT20,ST21} this problem was considered when $n=2$ for various subclasses of $\mathcal{S}$. In the present paper, we obtain the sharp bounds for $|A_2|-|A_1|$ and $|A_3|-|A_2|$ for functions belongs to some important subclasses of $\mathcal{S}$.

In this sequence, we have some subclasses of $\mathcal{S}$, which have been widely used by many authors in different prospective.

\subsection{The class $\mathcal{G}(\nu)$} In this paper, we also consider the class $\mathcal{G}(\nu)$. A locally univalent function $f\in \mathcal{A}$ is said to belong to $\mathcal{G}(\nu)$ for some $\nu>0$, if it satisfies the condition
\begin{equation*}
{\rm Re}\bigg(1+\cfrac{zf''(z)}{f'(z)}\bigg)<1+\cfrac{\nu}{2}, \quad z\in \mathbb{D}.
\end{equation*}
Ozaki \cite{Ozaki41} introduced the class $\mathcal{G}(1)=:\mathcal{G}$ and proved that functions in $\mathcal{G}$ are univalent in $\mathbb{D}$. Later Umezawa \cite{UME52} studied the class $\mathcal{G}$ and showed that this class contains the class of functions convex in one direction. Moreover, functions in $\mathcal{G}$ are proved to be starlike in $\mathbb{D}$ (see \cite{ponnusamy95},\cite{ponnusamy07}). Thus, the class $\mathcal{G}(\nu)$ is included in $\mathcal{S}^*$ whenever $\nu\in (0,1]$. It can be easily seen that functions in $\mathcal{G}(\nu)$ are not necessarily univalent in $\mathbb{D}$ if $\nu>1$. Recently, the radius of convexity for functions in the class $\mathcal{G}(\nu)$, $\nu>0$, is studied in \cite{Shankey20}.

\subsection{The class $\mathcal{F}(\lambda)$} For $-1/2< \lambda\le 1$, the class $\mathcal{F}(\lambda)$ defined  by 
\begin{equation*}
	\mathcal{F}(\lambda)=\left\{f\in \mathcal{A}:\,	{\rm Re}\left(1+\frac{zf''(z)}{f'(z)}\right)>\frac{1}{2} -\lambda ~\mbox{ for $z \in \mathbb{D}$}\right\}.
\end{equation*}
We note that clearly $\mathcal{F}(1/2)=:\mathcal{C}$ is the usual class of convex functions. Moreover, for $\lambda=1$, we obtain the class $\mathcal{F}(1)=:\mathcal{C}(-1/2)$ which considered by many researcher in the recent years. Also, functions in $\mathcal{C} (-1/2)$ are not necessarily starlike but are convex in some
direction. Other related results for $f\in \mathcal{C} (-1/2)$ were also obtained in \cite{AS18, PSY14}. Functions in $\mathcal{F}(\lambda)$  are close-to-convex for $1/2\le \lambda\le 1$ but $\mathcal{F}(\lambda)$ contains non-starlike functions for all $1/2< \lambda\le 1$ (see \cite{PRU76}). 
The class $\mathcal{F}(\lambda)$ was also considered for the restriction $1/2\le \lambda\le 1$, denote by $\mathcal{F}_0 (\lambda)$, and further extensively studied in the literature (see for instance \cite{ALT,KS20}).

\subsection{The class ${\mathcal C}_\gamma (\alpha)$} The family $\mathcal{C}_\gamma (\alpha )$ of  $\gamma$-convex functions of order $\alpha$ is defined by
$$
{\mathcal C}_\gamma (\alpha)=\left\{f\in \mathcal{A}:{\rm Re } \bigg( e^{-i\gamma }\left ( 1+\frac{zf''(z)}{f'(z)}\right )\bigg)>\alpha \cos \gamma\right\}
$$
where $0\le \alpha <1$ and $-\pi/2<\gamma<\pi/2$. We may set
${\mathcal C}_0(\alpha)=:{\mathcal C}(\alpha)$ which consists of the normalized
convex functions of order $\alpha$. A function in ${\mathcal C}_\gamma (0)$ need not be
univalent in $\mathbb{D}$ for general values of $\gamma $
$(|\gamma|<\pi/2)$. For example, the function
$f(z)=i(1-z)^i-i$ is known to belong to
${\mathcal C}_{\pi/4}\backslash {\mathcal S}$. Robertson
\cite{Robertson-69} showed  that $f\in\mathcal{C}_{\gamma}$ is
univalent if $0<\cos \gamma\leq 0.2315\cdots$. Finally, Pfaltzgraff \cite{Pfaltzgraff} has shown that
$f\in\mathcal{C}_{\gamma}$ is univalent whenever $0<\cos \gamma\leq
1/2$. This settles the improvement of range of $\gamma$ for which $f\in\mathcal{C}_{\gamma}$
is univalent. On the other hand,  in \cite{SinghChic-77} it was also shown that
functions in ${\mathcal C}_\gamma$ which satisfy $f''(0)=0$ are
univalent for all real values of $\gamma$ with $|\gamma|<\pi /2$.

Let $\mathcal{P}$ denote the class of all analytic functions $p$ having positive real part in $\mathbb{D}$, with the form
\begin{equation}\label{p}
p(z)=1+c_1z+c_2z^2+\cdots.
\end{equation}
A member of $\mathcal{P}$ is called a {\em Carath\'{e}odory function}. It is known that $|c_n|\le 2$ for a function $p\in \mathcal{P}$ and for all $n\ge 1$ (see \cite{Dur83}).\\[2mm]

To prove our results, we need the following lemma.

\begin{lemma}\cite{SimT20}\label{lemma2}
Let $B_1$, $B_2$, and $B_3$ be numbers such that $B_1>0$, $B_2\in \mathbb{C},$ and $B_3\in \mathbb{R}$. Let $p\in \mathcal{P}$ be of the form \eqref{p}. Define $\Psi_+(c_1,c_2)$ and $\Psi_-(c_1,c_2)$ by
$$
\Psi_+(c_1,c_2)=|B_2c_1^2+B_3c_2|-|B_1c_1|,
$$
and 
$$
\Psi_-(c_1,c_2)=-\Psi_+(c_1,c_2).
$$
Then
\begin{equation}\label{B+}
\Psi_+(c_1,c_2)\le\left \{
	\begin{array}{ll}
		|4B_2+2B_3|-2B_1, & {\mbox{ if }} |2B_2+B_3|\ge |B_3|+B_1,
		\\[5mm]
		2|B_3|, & {\mbox{ otherwise}},
	\end{array}
	\right.
\end{equation}
and
\begin{equation}\label{B-}
	\Psi_-(c_1,c_2)\le\left \{
	\begin{array}{ll}
		2B_1-B_4, & {\mbox{ if }} B_1\ge B_4+2|B_3|,
		\\[5mm]
		2B_1\sqrt{\cfrac{2|B_3|}{B_4+2|B_3|}}, & {\mbox{ if }} B_1^2\le 2|B_3|(B_4+2|B_3|),\\[5mm]
		2|B_3|+\cfrac{B_1^2}{B_4+2|B_3|}, &{\mbox{ otherwise}},
	\end{array}
	\right.
\end{equation}
where $B_4=|4B_2+2B_3|$. All inequalities in \eqref{B+} and \eqref{B-} are sharp.
\end{lemma}

Our main aim of this paper is to estimate the sharp bounds of $|A_2|-|A_1|$ and $|A_3|-|A_2|$ for functions $f$ belong to $\mathcal{G}(\nu)$, $\mathcal{F}_0(\lambda)$, and $\mathcal{C}_{\gamma}(\alpha)$. The organization of this paper is as follows: Section \ref{sect2} is devoted to the statements of main results. The proof of main results are
given in Section \ref{proofs}.

\section{Main results}\label{sect2}

We now state our first main result which provides sharp bounds for $|A_2|-|A_1|$ when $f$ belongs the class $\mathcal{G}(\nu)$.
\begin{theorem} \label{a2a1}
	Let $0< \nu\le 1$. For every $f\in \mathcal{G}(\nu)$ of the form \eqref{S}, we have
	$$
-1\le|A_2|-|A_1|\le \cfrac{2\lambda-1}{2}.
$$
Both inequalities are sharp.
\end{theorem}

\begin{theorem}\label{a3a2}
 Let $0< \nu \le 1$. For every $f\in \mathcal{G}(\nu)$ of the form \eqref{S}, we have
 $$
 |A_3|-|A_2|\le \cfrac{\nu}{6}
 $$
 and
\begin{equation}\label{gA3A2}
|A_3|-|A_2|\ge  \left \{
\begin{array}{ll}
	-\cfrac{\nu(8\nu+17)}{48(\nu+1)}, & {\mbox{ for }} 0<\nu\le 1/8,
	\\[5mm]
	\cfrac{-\nu}{2\sqrt{2(\nu+1)}}, & {\mbox{ for }} 1/8\le \nu\le 1.
\end{array}
\right.
\end{equation}
The inequalities are sharp.
\end{theorem}

Next, we obtain the sharp bounds for $|A_2|-|A_1|$ and $|A_3|-|A_2|$ when the functions $f$ are in $\mathcal{F}_0(\lambda)$.
\begin{theorem}\label{ozkia1a2}
Let $1/2\le \lambda\le 1$. For every $f\in \mathcal{F}_0(\lambda)$ be of the form \eqref{S}, we have
$$
-1\le|A_2|-|A_1|\le \cfrac{2\lambda-1}{2}.
$$
The inequalities are sharp.
\end{theorem}

\begin{theorem}\label{ozkia2a3}
Let $1/2\le \lambda\le 1$. For every $f\in \mathcal{F}_{0}(\lambda)$ of the form \eqref{S}, we have
\begin{equation}\label{A2A3}
-\cfrac{\sqrt{2\lambda+1}}{2\sqrt{2}}\le|A_3|-|A_2|\le \left \{
\begin{array}{ll}
	\cfrac{2\lambda+1}{6}, & {\mbox{ for }} 1/2\le \lambda\le 3/4,
	\\[5mm]
	\cfrac{(2\lambda+1)(2\lambda-1)}{3}, & {\mbox{ for }} 3/4\le \lambda\le 1.
\end{array}
\right.
\end{equation}
The inequalities are sharp.
\end{theorem}

In the next theorem, we will discuss about the sharp bounds for $|A_2|-|A_1|$ and $|A_3|-|A_2|$ when the functions $f$ are $\gamma $ -convex of order $\alpha$.
\begin{theorem}\label{convexa1a2}
	Let $-\pi/2< \gamma< \pi/2$ and $0\le \alpha <1$. For every $f\in \mathcal{C}_{\gamma}(\alpha)$ be of the form \eqref{S}, we have
	$$
	-1\le|A_2|-|A_1|\le (1-\alpha)\cos \gamma-1.
	$$
Both inequalities are sharp.	
\end{theorem}

\begin{theorem}\label{convexa2a3}
	Let $-\pi/2< \gamma< \pi/2$ and $0\le \alpha <1$. For every $f\in \mathcal{C}_{\gamma}(\alpha)$ of the form \eqref{S}, we have
	\begin{equation}\label{eqc}
	|A_3|-|A_2|\le \cfrac{(1-\alpha)\cos \gamma}{3}
	\end{equation}
and
\begin{equation}\label{eqco}
|A_3|-|A_2|\ge  \left \{
\begin{array}{ll}
	-\cfrac{(1-\alpha)\cos \gamma}{\sqrt{|\tau|+1}}, & {\mbox{ for }} |\tau|\ge 5/4,
	\\[5mm]
-(1-\alpha)\cos \gamma	\cfrac{13+4|\tau|}{12(|\tau|+1)}, & {\mbox{ for }} |\tau|\le 5/4,
\end{array}
\right.
\end{equation}
where $\tau:=4(1-\alpha)\mu-1$.
Both inequalities are sharp.
\end{theorem}
If we put $\alpha=0$ and $\gamma=0$ in Theorem \ref{convexa2a3}, then we obtain the following result for the class of convex functions:
\begin{corollary}
	For every $f\in \mathcal{C}$ of the form \eqref{S}, we have
	\begin{equation}\label{C}
	\cfrac{1}{2}\le	|A_3|-|A_2|\le \cfrac{1}{3}.
	\end{equation}
Both inequalities are sharp.
\end{corollary}
\section{Proof of the main results}\label{proofs}
This section is devoted to the detailed discussion on our proof of the main results.
\subsection{Proof of Theorem \ref{a2a1}}
 Let $f\in \mathcal{G}(\nu)$. Then there exists a function $p(z)=1+c_1z+c_2z^2+\cdots\in\mathcal{P}$ satisfying 
 \begin{equation}\label{gp(z)}
 	p(z)=\cfrac{1}{\nu} \bigg(\nu- \cfrac{2zf''(z)}{f'(z)}\bigg).
 \end{equation}
After writing  $f$ and $p$ in the series form and by comparing the coefficients of $z$ and $z^2$ in the above equation, we obtain the relations
\begin{equation}\label{eqqga2a3}
	a_2=-\cfrac{\nu c_1}{4} \mbox{ and } a_3=\cfrac{\nu^2 c_1^2-2\nu c_2}{24}.
\end{equation}
Thus from equation \eqref{eqqga2a3} and \eqref{A2}, we have
$$
|A_2|-|A_1|=\cfrac{\nu|c_1|}{4}-1\le \cfrac{\nu-2}{2},
$$
	where the last inequality is obtained by using $|c_n|\le 2$ for $n\ge 1$. 
	For the equality, let us consider the function $g_1 \in \mathcal{G}(\nu)$ satisfying \eqref{gp(z)} with $p_1(z)=(1+z)/(1-z)$.
	Then we have 
	\begin{equation*}
	g_{1}(z)=\cfrac{(1+z)^{1+\nu}-1}{\nu+1},\quad z\in \mathbb{D},
	\end{equation*}
	 for which $A_2=-\nu/2$ and $A_1=1$.
    On the other hand, 
	$$
	|A_1|-|A_2|=1-\cfrac{\nu|c_1|}{4}\le 1.
	$$
	It is easy to see that equality holds when $g_2 \in \mathcal{G}(\nu)$ defined by \eqref{gp(z)} with 
	$p_2(z)=(1+z^2)/(1-z^2)$.
 In this case
 \begin{equation*}
 g_2(z)=\int_{0}^{z}(1-t^2)^{\nu/2}dt=z+\cfrac{\nu}{6}z^3+\cdots,\quad z\in \mathbb{D}
 \end{equation*}
and
\begin{equation}\label{g2}
	g_2^{-1}(w)=w-\cfrac{\nu}{6}w^3+\cdots,\quad w\in \mathbb{D}_{r_0}.
\end{equation}
	This completes the proof.\hfill{$\Box$}
\subsection{Proof of Theorem \ref{a3a2}} Let $f\in \mathcal{G}(\nu)$. Then from equation \eqref{A2} and \eqref{eqqga2a3} we obtain
	\begin{align}\label{eqB}
	|A_3|-|A_2|
	&=\cfrac{\nu}{12}\bigg(|c_2+\nu c_1^2|-3| c_1|\bigg)\nonumber\\
	&=\cfrac{\nu}{12}\bigg(|B_3c_2+B_2c_1^2|-|B_1 c_1|\bigg),
		\end{align}
	where 
	$$
	B_1:=3, B_2:=\nu, \mbox{ and } B_3:=1.
	$$ 
	As $|A_3|-|A_2|$ is invariant under rotation, to simplify the calculation we assume that $c_1=c\in [0,2]$. Therefore, we can apply Lemma \ref{lemma2}. A simple calculation shows that, when $0<\nu\le 1$, the first condition $|2B_2+B_3|\ge |B_3|+B_1$ for $\Psi_+(x,c)$ is not satisfied. Hence it follows from Lemma \ref{lemma2} and the equation \eqref{eqB} that
	$$
	|A_3|-|A_2|\le \cfrac{\nu}{12}(2B_3)= \cfrac{\nu}{6}.
	$$

 Here equality holds for $g_2^{-1}$ given by \eqref{g2} in which the coefficient of $w^2$ is $0$ and $w^3$ is $-\nu/6$. Thus, the right-hand equality of the theorem has been proved.
 
 We now proceed to prove the left-hand side inequality. By checking the conditions for $\Psi_-(x,c)$ in Lemma \ref{lemma2}, we conclude that $B_1^2\le 2|B_3|(B_4+2|B_3|)$ holds but $B_1\ge B_4+2|B_3|$ does not hold for $\nu \ge 1/8$.
 Thus, Lemma \ref{lemma2} together with equation \eqref{eqB} leads to desired inequality \eqref{gA3A2}.

We now show that the inequalities in \eqref{gA3A2} are sharp by constructing extreme functions for both cases.
For the case $1/8\le \nu\le 1$, we consider a function $g_3$ satisfying \eqref{gp(z)} with $p_3\in \mathcal{P}$ defined by 
$$
p_3(z)=\cfrac{1-z^2}{1-2sz+z^2}=1+2sz+(4s^2-2)z^2+(8s^3-6s)z^3+\cdots, 
$$
where $s=1/\sqrt{2(\nu+1)}$. Then it is easy to see that the coefficients of $p_3$ are given by $c_1=2s$ and $c_2=-2\nu/(\nu+1)$. From \eqref{eqB}, we obtain $|A_3|-|A_2|=-\nu/(2\sqrt{2\nu+2})$ and so the inequality \eqref{gA3A2} in this case is sharp for $g_3\in \mathcal{G}(\nu)$.

In the similar way we can see that the inequality \eqref{gA3A2} in the case $0<\nu\le 1/8$ is sharp for $g_4\in \mathcal{G}(\nu)$ defined by \eqref{gp(z)}, where
$$
p_4(z)=\cfrac{1-z^2}{1-2rz+z^2},
$$
where $r=3/(4\nu+4)$. Then it is easy to see that the coefficients of $p_3$ are given by $c_1=2r$ and $c_2=(1-16\nu-8\nu^2)/(2\nu+2)^2$. From \eqref{eqB}, we obtain $|A_3|-|A_2|=\nu(8\nu+17)/48(\nu+1)$, which completes the proof of Theorem \ref{a3a2}.

\hfill{$\Box$}
\subsection{Proof of Theorem \ref{ozkia1a2}}
	Suppose $f\in \mathcal{F}_0(\lambda)$. Then from definition we can write
	\begin{equation}\label{pf}
	1+\cfrac{zf''(z)}{f'(z)}=\bigg(\cfrac{1}{2}+\lambda\bigg)p(z)+\cfrac{1}{2}-\lambda, \quad z\in \mathbb{D}.
	\end{equation}
	By using the Taylor series representations of the functions $f$ and $p$, and comparing the coefficients of $z^n\, (n=1,2)$ both the sides, we obtain
\begin{equation}\label{eqa2c2}
	a_2=\cfrac{(1+2\lambda)c_1}{4}\text{ and } a_3=\cfrac{(1+2\lambda)(2c_2+(1+2\lambda)c_1^2)}{24}.
\end{equation}
By using \eqref{A2} together with \eqref{eqa2c2} we can write
$$
|A_2|-|A_1|=\cfrac{(2\lambda+1)|c_1|}{4}-1\le \cfrac{2\lambda-1}{2}.
$$
The last inequality holds since $|c_1|\le 2$. In order to show that the inequality is sharp, first we consider the function $f_1\in \mathcal{F}_0(\lambda)$ defined by 	\begin{equation}\label{f1}
	f_1(z)=\cfrac{(1-z)^{-2\lambda}-1}{2\lambda}, \quad z\in \mathbb{D}.
\end{equation}
Hence, $f_1^{-1}$ is given by
$$
f_1^{-1}(w)=w-\cfrac{(1+2\lambda)}{2}w^2+\cfrac{(2\lambda+1)(4\lambda+1)}{6}w^3+\cdots
$$
for $w\in \mathbb{D}_{r_0}$. Thus, $|A_2|-|A_1|=(2\lambda-1)/2$, which shows that right hand side inequality is sharp.

 Secondly, we estimate the upper bound for $|A_1|-|A_2|=1-(2\lambda+1)| c_1|/2\le 1$. For the sharpness, let us consider the function $f_2$ given by
 \begin{equation}\label{f2}
 	f_2(z)=z+\cfrac{(2\lambda+1)}{6}z^3+\cdots
 \end{equation} and the corresponding inverse function is of the form
 $$
  f_2^{-1}(w)=w-\cfrac{(2\lambda+1)}{6}w^3+\cdots \quad w\in \mathbb{D}_{r_0}.
 $$
 Hence $|A_2|-|A_1|=-1$. This completes the proof.
\hfill{$\Box$}

\subsection{Proof of Theorem \ref{ozkia2a3}} Let $f\in \mathcal{F}_{0}(\lambda)$. Then by means of equation \eqref{A2} and \eqref{eqa2c2}, we see that
\begin{align}\label{A3A2}
|A_3|-|A_2|&=\cfrac{2\lambda+1}{24}\bigg(|(4\lambda+2)c_1^2-2c_2|-6|c_1|\bigg)\nonumber\\
&=\cfrac{2\lambda+1}{24}(|B_2c_1^2+B_3c_2|-|B_1c_1|),
\end{align}
where 
$$
B_1:=6, B_2:=4\lambda+2\mbox{ and }B_3:=-2.
$$
We can see that the functional $|A_3|-|A_2|$ is rotationally invariant, so we assume $c_1=c\in [0,1]$. Thus, we can apply Lemma \ref{lemma2} and by checking the conditions for the bound $\Psi_+(x,c)$ we obtain 
\begin{equation*}
	\Psi_+(c_1,c_2)\le\left \{
	\begin{array}{ll}
		|4B_2+2B_3|-2B_1, & {\mbox{ for }} 3/4\le \lambda\le 1,
		\\[5mm]
		2|B_3|, & {\mbox{ for }} 1/2\le \lambda\le 3/4.
	\end{array}
	\right.
\end{equation*}
Therefore, Lemma \eqref{lemma2} and equation \eqref{A3A2} yields
\begin{equation*}
	|A_3|-|A_2|\le\left \{
	\begin{array}{ll}
		\cfrac{2\lambda+1}{6}, & {\mbox{ for }} 1/2\le \lambda\le 3/4,
		\\[5mm]
		\cfrac{(2\lambda+1)(2\lambda-1)}{3}, & {\mbox{ for }} 3/4\le \lambda\le 1.
	\end{array}
	\right.
\end{equation*}

We next find the lower bound of $|A_3|-|A_2|$. We can apply Lemma \ref{lemma2} for $\Psi(c_1,c_2)$ and we obtain that the condition $B_1^2\le 2|B_3|(B_4+2|B_3|)$ satisfies for $\lambda\in [1/2,1]$. Using Lemma \ref{lemma2} and equation \eqref{A3A2}, we obtain 
\begin{equation}\label{-A2A3}
|A_3|-|A_2|\ge- 2B_1\cfrac{2\lambda+1}{24} \sqrt{\cfrac{2|B_3|}{B_4+2|B_3|}}=-\cfrac{\sqrt{2\lambda+1}}{2\sqrt{2}}.
\end{equation}

We now show that inequalities are sharp. When $1/2\le \lambda\le 3/4$, equality holds in \eqref{A2A3} for $f_2$ given by \eqref{f2} and for $3/4\le \lambda\le 1$ equality holds in \eqref{A2A3} for $f_1$ given by \eqref{f1}. For the left-hand side equality, let us consider the function $f_3 \in \mathcal{F}_0(\lambda)$ satisfying \eqref{pf} with 
$$
p_3(z)=\cfrac{1+2tz+z^2}{1-z^2},
$$
where $t=1/\sqrt{4\lambda+2}$, for which 
 $$
 f_3^{-1}(w)=w-\cfrac{\sqrt{1+2\lambda}}{2\sqrt{2}}\,w^2+0\cdot w^3+\cdots.
 $$
This completes the proof.
 \hfill{$\Box$}	

\subsection{Proof of Theorem \ref{convexa1a2}}
Suppose $f\in \mathcal{C}_{\gamma}(\alpha)$. Then there exists a function $p\in \mathcal{P}$ such that
\begin{equation}\label{convexp(z)}
p(z)=\cfrac{1}{1-\alpha}\bigg\{ \cfrac{1}{\cos \gamma}\bigg(e^{-i\gamma} \bigg(1+\cfrac{zf''(z)}{f'(z)}\bigg)+i\sin \gamma\bigg)-\alpha\bigg\}.
\end{equation}
Equating the coefficients of $z^n$ on both the sides of \eqref{convexp(z)} for $n=1,2$, we obtain
\begin{equation}\label{eqa2a3}
	a_1=1,\,2a_2=(1-\alpha)\mu c_1 \mbox{  and  } 6a_3=(1-\alpha)^2\mu^2 c_1^2+(1-\alpha)\mu c_2,
\end{equation}
where $\mu = e^{i \gamma}\cos \gamma$. Thus from \eqref{A2} and \eqref{eqa2a3} we obtain
$$
|A_2|-|A_1|=\cfrac{(1-\alpha) |\mu c_1|}{2}-1\le (1-\alpha)\cos \gamma-1.
$$
The last inequality holds since $|c_1|\le 2$. It is easy to see that the equality holds for the function
$$
	h_{1}(z)=\cfrac{1}{2(1-\alpha)\mu-1}\bigg(\cfrac{1}{(1-z)^{2(1-\alpha)\mu-1}}-1\bigg)=z+(1-\alpha)\mu z^2+\cdots,\quad z\in \mathbb{D}
	$$
	and corresponding inverse function is
$$
h_{1}^{-1}(w)=w-(1-\alpha)\mu w^2+\cdots,\quad z\in \mathbb{D}_{r_0}
$$
	
Next, we compute the upper bound for $|A_1|-|A_2|=1-(1-\alpha)|\mu c_1|/2\le 1$. Now we shall easily see that the inequality is sharp for the function
$$
h_2(z)=z+\cfrac{(1-\alpha)\mu}{3}\,z^3+\cfrac{(1-\alpha)\mu((1-\alpha)\mu+1)}{10}\,z^5+\cdots
$$
and
\begin{equation}\label{h2}
h_2^{-1}(w)=w-\cfrac{(1-\alpha)\mu}{3}\,w^3+\cdots.
\end{equation}
This completes the proof.
 \hfill{$\Box$}	

\subsection{Proof of Theorem \ref{convexa2a3}}
Let $f\in \mathcal{C}_{\gamma}(\alpha)$. Then from \eqref{eqa2a3} and \eqref{A2}, we obtain 
\begin{align}\label{cB2}
	|A_3|-|A_2|&=\cfrac{(1-\alpha)\cos \gamma}{6}\bigg(|2(1-\alpha)\mu c_1^2-c_2|-3|c_1|\bigg)\nonumber\\
	&=\cfrac{(1-\alpha)\cos \gamma}{6}\bigg(|B_2 c_1^2+B_3c_2|-|B_1c_1|\bigg),
\end{align}
where
$$
B_1:=3,\,B_2:=2(1-\alpha)\mu, \mbox{ and }B_3=-1.
$$
Since $|A_3|-|A_2|$ is rotationally invariant, we may assume that $c_1=c\in [0,2]$. Therefore, we can apply Lemma \ref{lemma2}. A simple calculation shows that the first condition $|2B_2+B_3|\ge |B_3|+B_1$ for $\Psi_+(x,c)$ is not satisfied since $|4(1-\alpha)\mu-1|\le 3$. 
By using Lemma \ref{lemma2} and the equation \eqref{cB2}, we obtain that 
$$
|A_3|-|A_2|\le \cfrac{(1-\alpha)\cos \gamma}{3}.
$$
This proves the inequality \eqref{eqc}.

We next prove the lower bound in \eqref{eqco} by checking the condition of Lemma \ref{lemma2} for $\Psi_-(c_1,c_2)$. Note that the inequality $B_1\ge B_4+2|B_3|$ does not hold as $|4(1-\alpha)\mu-1|\le 3$
and $ B_1^2\le 2|B_3|(B_4+2|B_3|)$ hold when $|4(1-\alpha)\mu-1|\ge 5/4$. Thus, we can apply Lemma \ref{lemma2}, we obtain
$$
\Psi_-(c_1,c_2)\le\left \{
\begin{array}{ll}
	2B_1\sqrt{\cfrac{2|B_3|}{B_4+2|B_3|}}, & {\mbox{ if }} |4(1-\alpha)\mu-1|\ge 5/4,\\[5mm]
	2|B_3|+\cfrac{B_1^2}{B_4+2|B_3|}, &{\mbox{ if }} |4(1-\alpha)\mu-1|\le 5/4.
\end{array}
\right.
$$
By Substituting the above inequality in \eqref{cB2} we obtain the required inequality \eqref{eqco}.

 It is easy to see that equality holds in \eqref{eqc} when $h_2^{-1}$ is defined by \eqref{h2}.

For $|\tau|\ge 5/4$, we construct a function $h_3$ defined by \eqref{convexp(z)} with
\begin{equation}\label{eqcp(z)}
p(z)=\cfrac{1+q_1(q_2+1)z+q_2z^2}{1+q_1(q_2-1)-q_2z^2},
\end{equation}
where
$$
q_1=\cfrac{1}{\sqrt{|\tau|+1}} \mbox{ and } q_2=e^{i\arg \tau}
$$
with $\tau=4(1-\alpha)\mu-1$. Then
  $A_2=-(1-\alpha)\mu/\sqrt{|\tau|+1}$ and $A_3=0$ which gives the equality in \eqref{eqco}.
  
 For $|\tau|\le 5/4$, consider the function $h_4$ given by \eqref{eqcp(z)} where $p(z)$ is of the form \eqref{eqcp(z)} with 
 $$
 q_1=\cfrac{3}{2(|\tau|+1)} \mbox{ and } q_2=e^{i\arg \tau}.
 $$
 Then we have
 $$
 A_2=-\cfrac{3(1-\alpha)\mu}{2(|\tau|+1)} \mbox{ and } A_3=\cfrac{(1-\alpha)\mu\tau(-4|\tau|+5)}{12|\tau|(|\tau|+1)}.
 $$
 This completes the proof of this theorem.
 \hfill{$\Box$}	

\bigskip
\noindent
{\bf Acknowledgement.}
The first author thank SERB-CRG and the second author thank IIT Bhubaneswar for providing Institute Post Doctoral fellowship.

\end{document}